\documentclass[11pt,a4paper]{article}
\usepackage[T1]{fontenc}
\usepackage[utf8]{inputenc}
\usepackage[english]{babel}
\usepackage{amsmath,amssymb,amsthm,mathtools}
\usepackage[margin=27mm]{geometry}
\usepackage[hidelinks]{hyperref}
\newtheorem{theorem}{Theorem}[section]
\newtheorem{proposition}[theorem]{Proposition}
\newtheorem{lemma}[theorem]{Lemma}
\newtheorem{corollary}[theorem]{Corollary}
\theoremstyle{remark}
\newtheorem{remark}[theorem]{Remark}
\DeclareMathOperator{\dist}{dist}
\newcommand{\HH}{\mathcal H}
\newcommand{\CC}{\mathcal C}
\newcommand{\TT}{\mathcal T}
\newcommand{\EE}{\mathcal E}
\title{Multifractal analysis of Takagi-type functions built from beta-expansion cylinders}
\author{S. Nicolay\footnote{Universit\'e de Li\`ege, D\'epartement de math\'ematique, zone Polytech 1, 12 all\'ee de la D\'ecouverte, B\^at. B37, B-4000 Li\`ege. Orcid ID: 0000-0003-0549-0566. S.Nicolay@uliege.be}}
\date{\today}
\begin{document}
\maketitle
\begin{abstract}
We study a family of Takagi functions obtained by placing a normalized tent on each beta-expansion cylinder, with amplitude $\beta^{-\alpha n}$ at order $n$, where $\beta>1$ and $0<\alpha\le1$. We characterize their pointwise H\"older regularity in terms of cylinder lengths and distances from the point under consideration. The characterization holds more generally for nested interval partitions with exponentially decreasing mesh, without requiring interval lengths at a given order to be comparable. For beta-expansion cylinders, we compute the full multifractal spectrum whenever the remainders of the expansion of $1$ have a finite upper exponential decay rate. The functions are monofractal when this rate vanishes, as it does for Parry, Pisot and Salem bases. When it is positive, the spectrum is affine on a nondegenerate interval. The bases defined by $\sum_{k\ge0}\beta_\ell^{-\ell^k}=1$, for integers $\ell\ge2$, give explicit non-monofractal examples. We also prove that every function in the family is nowhere differentiable.
\end{abstract}

\noindent Keywords: Takagi-type functions; beta-expansions; pointwise H\"older regularity; multifractal spectrum; Hausdorff dimension.\\
\noindent 2020 MSC: 26A16; 28A80, 37E05, 11A63.

\section*{Introduction}
The Takagi function is a classical example of continuous nowhere differentiable function, defined by a series of piecewise affine functions. Writing $\phi(t)=\dist(t,\mathbb Z)$, we consider the family
\[
 \sum_{n\ge0}m^{-\alpha n}\phi(m^n x),
 \]
where $m\ge2$ is an integer and $\alpha\in (0,1]$. At order $n$, the function $\phi(m^n x)$ consists of tents of equal height supported on intervals of length $m^{-n}$. Such constructions allow one to study pointwise regularity and the Hausdorff dimensions of its level sets. Related questions arise for other classical irregular functions: Jaffard and Nicolay \cite{JN} determine the pointwise regularity of Lebesgue's space-filling function and related examples, including both monofractal and multifractal functions. Generalizations of Takagi's construction include the self-affine functions studied by Allaart \cite{Allaart} and functions arising from random dynamics, whose multifractal spectra are determined by Jaerisch and Sumi \cite{JS}.

Here we retain the tent construction and replace the integer-base subdivisions by the cylinders of expansions in a real base $\beta>1$. A cylinder of order $n$ is the interval of points with prescribed first $n$ digits. Let $\CC_n(\beta)$ denote the family of closed cylinders and let $\tau_I$ be the tent supported on $I$, vanishing at its endpoints and taking the value $1/2$ at its midpoint. We study
\[
 \TT_{\beta,\alpha}(x)
 =\sum_{n\ge0}\beta^{-\alpha n}
       \sum_{I\in\CC_n(\beta)}\tau_I(x),
\]
with  $\alpha\in (0,1]$. For integer $\beta$, this is the series above. For a real base, cylinders of the same order can have very different lengths. Their tents still have the same height, so the shortest cylinders carry the largest slopes. We determine how these variations in length affect the regularity of the sum.

The length of the cylinder containing $x$ does not by itself determine the pointwise exponent. Nearby cylinders also matter, through the quantity
\[
 \delta_n(x)=\min_{I\in\CC_n(\beta)}
                  \bigl(|I|+\dist(x,I)\bigr).
\]
We obtain an exact formula for the exponent in terms of $\delta_n(x)$. The proof applies to finite nested partitions with exponentially decreasing mesh, without requiring interval lengths at a fixed order to be comparable (Theorem~\ref{thm:partitions}). The same geometric argument proves nowhere differentiability (Corollary~\ref{cor:nondifferentiability}).

For beta-expansions, cylinder lengths are governed by the remainders in the greedy expansion of the number $1$. Their upper exponential decay rate defines the parameter $\lambda(\beta)$ in Section~\ref{sec:spectra}. When this parameter is finite, we compute the full spectrum of $\TT_{\beta,\alpha}$ (Theorem~\ref{thm:general}). The case $\lambda(\beta)=0$ is monofractal and includes Parry bases, as well as Pisot and Salem bases by an algebraic estimate. If $\lambda(\beta)>0$, the spectrum is affine on a nondegenerate interval. The bases $\beta_\ell$ defined by
\[
 \sum_{k\ge0}\beta_\ell^{-\ell^k}=1
\]
for $\ell\ge2$, give explicit examples with $\lambda(\beta_\ell)=\ell-1$ (Section~\ref{sec:family}).

Suzuki \cite{Suzuki} introduces a different generalization of the Takagi function using the digits of real-base expansions; we compare the two
definitions in Section~\ref{sec:cylinders}. The spectrum theorems of Allaart and Jaerisch--Sumi do not directly yield our results. In the setting of positive horizontal contractions in Allaart's Theorem~1.4 \cite{Allaart}, a nondegenerate spectrum is strictly concave in case~(a), whereas case~(b), which allows an affine part, has no finite exponent below $1$. Neither case gives the affine spectrum of our non-monofractal examples. Under the open set condition, Theorem~1.3 of Jaerisch--Sumi \cite{JS} likewise gives a strictly concave spectrum in the nondegenerate case.

Wavelet criteria \cite{AJ,CN}, and criteria involving Schauder coefficients \cite{BS,BDK} provide a point of comparison. In these approaches, regularity is related to the size of coefficients, their supports, and their distance from the point under consideration. Here the function is defined directly as a series of tents with prescribed amplitudes, and our analysis uses those amplitudes together with the geometry of their supports. Theorem~3.4 of Bayraktar--Das--Kim \cite{BDK} is relevant in a different respect: it characterizes global regularity under a balance assumption on the partitions, an assumption not imposed here. The remarks following the definitions of the tents and of \(\delta_n(x)\) specify how these features enter our analysis.

Fan and Wang \cite{FW} determine the spectrum associated with the length of the cylinder containing a point. Their result is a natural antecedent, but does not account for the neighbouring cylinders occurring in $\delta_n(x)$. Section~\ref{sec:comparison} compares the two quantities and their level sets. Our spectrum calculation uses the distribution of full cylinders studied by Bugeaud and Wang \cite{BW} and the mass transference principle of Beresnevich and Velani \cite{BV}.

\section{The family and its pointwise regularity}\label{sec:partitions}
By a closed interval partition of $[0,1]$ we mean a finite family of closed intervals of positive length covering $[0,1]$, with pairwise disjoint interiors. Adjacent intervals may therefore share an endpoint. A partition refines another if each of its intervals is contained in an interval of the other partition.

For $I=[a,b]$, define the normalized tent
\[
\tau_I(x)=\begin{cases}
\displaystyle\min\{\frac{x-a}{b-a},\frac{b-x}{b-a}\}&\text{if }a\le x\le b\\
0&\text{if }x\notin[a,b]
\end{cases}.
\]
It is affine on each half of $I$, vanishes at the endpoints and equals $1/2$ at the midpoint. Its slopes are $1/|I|$ and $-1/|I|$.

Let $\rho>1$, $0<\alpha\le1$, and let $\mathcal P=(\mathcal P_n)_{n\ge0}$ be a sequence of closed interval partitions, each refining the preceding one, such that $\max_{I\in\mathcal P_n}|I|\le C_0\rho^{-n}$ for a constant $C_0>0$. Define
\begin{equation}\label{eq:def-general}
\TT_{\mathcal P,\rho,\alpha}(x)
=\sum_{n\ge0}\rho^{-\alpha n}\sum_{I\in\mathcal P_n}\tau_I(x).
\end{equation}
For each $x$ and $n$, at most one tent in the inner sum is nonzero. At partition endpoints, all tents of that order vanish. The inner sum
therefore lies between $0$ and $1/2$, so the series converges uniformly and defines a continuous function.

\begin{remark}
The form of \eqref{eq:def-general} suggests studying its multifractal properties through a Schauder basis. The tents used here must nevertheless be distinguished from the basis functions associated with successive subdivisions. The midpoint of an interval in $\mathcal P_n$ need not be a subdivision point at subsequent orders. Thus the amplitudes $\rho^{-\alpha n}$ cannot be identified with the Schauder coefficients of the sum without further justification. We prove the pointwise criterion directly from the values of the series.
\end{remark}

For a continuous function $f$, its pointwise H\"older exponent is
\[
h_f(x)=\sup\{s>0:\exists P,\ \deg P<s,\quad
f(y)-P(y-x)=O(|y-x|^s)\text{ as }y\to x\},
\]
where $P$ is a polynomial. Limits at the endpoints are relative to $[0,1]$, and the supremum is taken to be zero if no positive order is
admissible. The spectrum is
\[
D_f(u)=\dim_H\{x:h_f(x)=u\},
\]
where $\dim_H$ denotes the Hausdorff dimension, with the convention $\dim_H\varnothing=-\infty$. Since the exponent is a supremum, the estimate at the critical order need not hold.

For positive quantities, $A\lesssim B$ means $A\le CB$, and $A\asymp B$ means $cB\le A\le CB$, with $c,C>0$ independent of the varying indices and points. Constants may depend on fixed parameters; any additional dependence will be specified.

Set
\begin{equation}\label{eq:delta}
\delta_n(x)=\min_{I\in\mathcal P_n}(|I|+\dist(x,I)).
\end{equation}
The minimum is taken over all intervals in the partition, including those not containing $x$. If $I_n$ contains $x$, then $\delta_n(x)\le |I_n|$, but equality need not hold. Since $0<\delta_n(x)\le C_0\rho^{-n}$, the logarithmic denominators below are positive for all sufficiently large $n$.

\begin{remark}\label{sec:regularity-prior}
The quantity $|I|+\dist(x,I)$ also occurs in pointwise wavelet criteria. Proposition~1 of Ayache--Jaffard \cite{AJ}, p.~86, recalls that if
$f\in C^s(x)$, its $L^\infty$-normalized wavelet coefficients $c_{j,k}(f)$ are bounded in absolute value by a constant times $(2^{-j}+|x-2^{-j}k|)^s$ for indices near $x$. Here $f\in C^s(x)$ means that there is a polynomial $P$ of degree strictly less than $s$ such that $f(y)-P(y-x)=O(|y-x|^s)$ as $y\to x$. The same proposition gives a converse, with a logarithmic loss, under an additional decay assumption.
\end{remark}

\begin{lemma}
\label{lem:chord}
If $I=[a,b]\in\mathcal P_n$ and $c=(a+b)/2$, then
\begin{equation}\label{eq:chord}
 \TT_{\mathcal P,\rho,\alpha}(c)
 -\frac{\TT_{\mathcal P,\rho,\alpha}(a)+\TT_{\mathcal P,\rho,\alpha}(b)}2
 \ge\frac12\rho^{-\alpha n}.
\end{equation}
\end{lemma}
\begin{proof}
The contribution of
$g_n=\rho^{-\alpha n}\sum_{J\in\mathcal P_n}\tau_J$ to the left-hand side is $\rho^{-\alpha n}/2$: the tent $\tau_I$ equals $1/2$ at $c$, all other tents of order $n$ vanish there, and all vanish at $a,b$. For $k<n$, the interval $I$ is contained in an interval of $\mathcal P_k$, whose tent is concave on $I$ and hence has a nonnegative contribution. For $k>n$, both $a$ and $b$ are partition endpoints. All tents of order $k$ vanish there and are nonnegative at $c$. Uniform convergence allows these inequalities to be summed.
\end{proof}
\begin{theorem}\label{thm:partitions}
For every $x\in[0,1]$,
\begin{equation}\label{eq:holder-partitions}
h_{\TT_{\mathcal P,\rho,\alpha}}(x)
=H(x):=\liminf_{n\to\infty}
\frac{\alpha n\log\rho}{-\log\delta_n(x)}.
\end{equation}
In particular, $0\le h_{\TT_{\mathcal P,\rho,\alpha}}(x)\le\alpha$.
\end{theorem}
\begin{proof}
Since $\delta_n(x)\le C_0\rho^{-n}$, we have $H(x)\le\alpha\le1$. Suppose first that $h_{\TT_{\mathcal P,\rho,\alpha}}(x)>H(x)$ and choose $H(x)<s<\min\{h_{\TT_{\mathcal P,\rho,\alpha}}(x),2\}$ with $s\ne1$. There are a polynomial $P$ of degree less than $s$, a constant $C(\mathcal P,\rho,\alpha,x,s)>0$ and a radius $r(\mathcal P,\rho,\alpha,x,s)>0$ such that
\[
 |\TT_{\mathcal P,\rho,\alpha}(y)-P(y-x)|
 \le C(\mathcal P,\rho,\alpha,x,s)|y-x|^s
\]
for $y\in[0,1]$ with $|y-x|<r(\mathcal P,\rho,\alpha,x,s)$. This polynomial approximates the function near $x$. Since $s<2$, it is constant or affine and contributes zero to the deviation from the chord. Choose $I_n$ attaining the minimum in \eqref{eq:delta}. Every point of $I_n$ is at distance at most $\delta_n(x)$ from $x$. The polynomial approximation and Lemma~\ref{lem:chord} give, for sufficiently large $n$,
\[
 \tfrac12\rho^{-\alpha n}\le 2C(\mathcal P,\rho,\alpha,x,s)\delta_n(x)^s.
\]
Taking logarithms and the lower limit gives $H(x)\ge s$, a contradiction.

For the reverse inequality, suppose $H(x)>0$ and fix $0<s<H(x)$, so that $s<1$. By the definition of the lower limit, there is an integer $n_0(\mathcal P,\rho,\alpha,x,s)$ such that, for every $n\ge n_0(\mathcal P,\rho,\alpha,x,s)$,
\begin{equation}\label{eq:coeff}
 A_n:=\rho^{-\alpha n}\le\delta_n(x)^s
       \le\bigl(\dist(x,I)+|I|\bigr)^s
       \quad(I\in\mathcal P_n).
\end{equation}
Write $g_n=A_n\sum_{I\in\mathcal P_n}\tau_I$ and $t=|y-x|$. We claim that
\begin{equation}\label{eq:each}
 |g_n(y)-g_n(x)|\le (1+2^s)t^s
 \quad(n\ge n_0(\mathcal P,\rho,\alpha,x,s)).
\end{equation}
If $x,y$ belong to the same interval of length $L$, then $A_n\le L^s$ and
\[
|g_n(y)-g_n(x)|\le A_n\min\{1,t/L\}\le t^s.
\]
Otherwise, choose intervals $I_x,I_y$ containing $x,y$, of lengths $L_x,L_y$. Choices at endpoints do not matter because the tents vanish there. Each interval has an endpoint between $x$ and $y$, so
\[
 A_n\tau_{I_x}(x)\le L_x^s\min\{1,t/L_x\}\le t^s.
\]
Also $\dist(x,I_y)\le t$, hence $A_n\le(L_y+t)^s$. If $L_y\ge t$, then
\[
A_n\tau_{I_y}(y)\le 2^s L_y^s t/L_y\le2^s t^s.
\]
If $L_y<t$, its height is at most $A_n/2\le(2t)^s/2$. This proves \eqref{eq:each}.

The sum of the terms with $n<n_0(\mathcal P,\rho,\alpha,x,s)$ is Lipschitz, being a finite sum of piecewise affine functions. Its Lipschitz constant depends on $\mathcal P,\rho,\alpha,x,s$. For $0<t<1$, put $N(t,s,\rho,\alpha)=\lceil s\log(1/t)/(\alpha\log\rho)\rceil$; then $\rho^{-\alpha N(t,s,\rho,\alpha)}\le t^s$. Summing \eqref{eq:each} up to $N(t,s,\rho,\alpha)$ and using $\|g_n\|_\infty\le A_n/2$ for the tail gives
\[
 |\TT_{\mathcal P,\rho,\alpha}(y)-\TT_{\mathcal P,\rho,\alpha}(x)|
 \le C(\mathcal P,\rho,\alpha,x,s)t^s(1+\log(1/t)).
\]
For every $0<s'<s$, the right-hand side is $O(t^{s'})$ as $t\to0$. Thus $h_{\TT_{\mathcal P,\rho,\alpha}}(x)\ge s$; letting $s\uparrow H(x)$ proves the result. If $H(x)=0$, the lower bound is immediate.
\end{proof}

\begin{corollary}
\label{cor:nondifferentiability}
For $0<\alpha\le1$, the function $\TT_{\mathcal P,\rho,\alpha}$ has no finite derivative at any point of $(0,1)$ and no finite one-sided derivative at either endpoint of $[0,1]$.
\end{corollary}
\begin{proof}
Write $f=\TT_{\mathcal P,\rho,\alpha}$ and suppose that $f$ is differentiable at $x\in[0,1]$, with the one-sided interpretation at an endpoint. For each $n$, choose $I_n=[a_n,b_n]\in\mathcal P_n$ containing $x$, and set $c_n=(a_n+b_n)/2$. We have $|I_n|\to0$ and
\[
 f(y)=f(x)+f'(x)(y-x)+o(|y-x|).
\]
The points $a_n,b_n,c_n$ are at distance at most $|I_n|$ from $x$. The affine part cancels, giving
\[
 f(c_n)-\frac{f(a_n)+f(b_n)}2=o(|I_n|).
\]
On the other hand, Lemma~\ref{lem:chord} and $|I_n|\le C_0\rho^{-n}$ give
\[
 \frac{f(c_n)-\frac12(f(a_n)+f(b_n))}{|I_n|}
 \ge \frac{\rho^{-\alpha n}}{2|I_n|}
 \ge \frac{\rho^{(1-\alpha)n}}{2C_0}
 \ge \frac1{2C_0},
\]
a contradiction.
\end{proof}

\section{Beta-expansion cylinders}\label{sec:cylinders}
Let $\beta>1$. For $x\in[0,1)$, define
\[
T_\beta(x)=\beta x-\lfloor\beta x\rfloor,\qquad
 a_j(x)=\lfloor\beta T_\beta^{j-1}(x)\rfloor.
\]
The sequence $(a_j(x))_{j\ge1}$ consists of the digits of the greedy beta-expansion of $x$. The cylinder associated with a word $w=(w_1,\ldots,w_n)$ is
\[
 C(w)=\{x\in[0,1):a_j(x)=w_j\text{ for }1\le j\le n\}.
\]
We call $w$ admissible if $C(w)$ has positive length. Such a cylinder is a half-open interval $[a,b)$. Let $\CC_n(\beta)$ be the family of their closures $[a,b]$, and set $\CC_0(\beta)=\{[0,1]\}$. On every cylinder of order $n$, the map $T_\beta^n$ is affine with slope $\beta^n$ and image $[0,r)$ for some $r\in(0,1]$. The partitions $\CC_n(\beta)$ are therefore nested and have mesh at most $\beta^{-n}$.
A cylinder is full if its length is $\beta^{-n}$, or equivalently if its image under $T_\beta^n$ is $[0,1)$.

\begin{lemma}\label{lem:lower-digit}
Let $w=(w_1,\ldots,w_n)$ be admissible, with $n\ge1$ and $w_n\ge1$. The cylinder associated with $(w_1,\ldots,w_{n-1},w_n-1)$ is full and has right endpoint $\sum_{i=1}^{n}w_i\beta^{-i}$.
\end{lemma}
\begin{proof}
The image of the parent cylinder under $T_\beta^{n-1}$ is $[0,r)$. Since $C(w)$ has positive length, $\beta r>w_n$. Thus $[(w_n-1)/\beta,w_n/\beta)$ is contained in $[0,r)$ and is a full branch of $T_\beta$, with image $[0,1)$. The corresponding cylinder is full, and its right endpoint is $\sum_{i=1}^{n-1}w_i\beta^{-i}+w_n\beta^{-n}$.
\end{proof}

Apply the general construction with $\mathcal P_n=\CC_n(\beta)$ and $\rho=\beta$, and write
\begin{equation}\label{eq:def}
\TT_{\beta,\alpha}(x)=\sum_{n\ge0}\beta^{-\alpha n}
                  \sum_{I\in\CC_n(\beta)}\tau_I(x).
\end{equation}
For every integer $m\ge2$, this gives exactly
\[
\TT_{m,\alpha}(x)=\sum_{n\ge0}m^{-\alpha n}\phi(m^nx),
\]
with $\phi(t)=\dist(t,\mathbb Z)$. The choice $m=2$, $\alpha=1$ recovers the classical Takagi function. For noninteger bases, the cylinder construction generally differs from the series $\sum_n\beta^{-\alpha n}\phi(\beta^nx)$.

\begin{remark}
\label{rem:suzuki}
For $1<\beta<2$, formula~(2.2) of Suzuki \cite{Suzuki} gives $G_\beta(1/\beta)=1/\beta$. At $1/\beta$, all tents of positive order in our construction vanish, whereas the tent of order zero equals $1-1/\beta$. Hence
\[
\TT_{\beta,\alpha}(1/\beta)=1-1/\beta.
\]
In fact, there are no constants $a,b$ such that $\TT_{\beta,\alpha}=aG_\beta+b$. At $x=0$, both functions vanish, so such an identity would require $b=0$; evaluating at $1/\beta$ would then give $a=\beta-1$. For $m\ge1$, the greedy expansion of $\beta^{-m}$ has a single nonzero digit, equal to $1$ at position $m$. Suzuki's formula~(2.2) therefore gives $G_\beta(\beta^{-m})=m\beta^{-m}$. On the other hand, $\beta^{-m}$ belongs to $[0,\beta^{-n}]$ for $n<m$. Setting $t_j=\min\{\beta^{-j},1-\beta^{-j}\}$, we obtain
\[
\TT_{\beta,\alpha}(\beta^{-m})
 =\beta^{-\alpha m}\sum_{j=1}^{m}\beta^{\alpha j}t_j.
\]
If $0<\alpha<1$, the sum tends to a finite positive constant, since $t_j=\beta^{-j}$ for all sufficiently large $j$. Consequently $\TT_{\beta,\alpha}(\beta^{-m})/G_\beta(\beta^{-m})\to+\infty$. If $\alpha=1$, the sum is $m+O(1)$ and the quotient tends to $1$. In neither case can the quotient be constantly equal to $\beta-1$.
\end{remark}

\begin{corollary}\label{thm:holder}
For every $x\in[0,1]$,
\begin{equation}\label{eq:holder}
h_{\TT_{\beta,\alpha}}(x)=\liminf_{n\to\infty}
\frac{\alpha n\log\beta}{-\log\delta_n(x)}.
\end{equation}
\end{corollary}
\begin{proof}
Apply Theorem~\ref{thm:partitions} to the partitions $\CC_n(\beta)$, with $C_0=1$.
\end{proof}

To evaluate a partial sum at $x$, only the cylinder containing $x$ at each order is needed. The following recurrence computes the normalized position of $x$ within that cylinder without constructing the entire partition. It applies to every $\beta>1$, including the bases $\beta_\ell$ of Section~\ref{sec:family}. For $x\in[0,1)$, set $y_0=x$, $r_0=1$, and
\begin{equation}\label{eq:algorithm}
 d_{n+1}=\lfloor\beta y_n\rfloor,\qquad
 y_{n+1}=\beta y_n-d_{n+1},\qquad
 r_{n+1}=\min\{1,\beta r_n-d_{n+1}\}.
\end{equation}
Then
\begin{equation}\label{eq:explicit}
 \TT_{\beta,\alpha}(x)=
 \sum_{n\ge 0}\beta^{-\alpha n}
      \min\left\{\frac{y_n}{r_n},1-\frac{y_n}{r_n}\right\}.
\end{equation}
At the remaining endpoint, $\TT_{\beta,\alpha}(1)=0$. Indeed, $y_n=T_\beta^n(x)$ and $[0,r_n)$ is the image of the cylinder containing $x$ at order $n$. Since the restriction of $T_\beta^n$ is affine with slope $\beta^n$, the quotient $y_n/r_n$ is the normalized position of $x$ in the cylinder. Subdividing $[0,r_n)$ according to the digit $d_{n+1}$ gives the last relation in \eqref{eq:algorithm}.

If $S_N$ is the sum over orders $0,\ldots,N-1$, then
\begin{equation}\label{eq:tailnum}
 0\le\TT_{\beta,\alpha}(x)-S_N(x)
 \le \frac{\beta^{-\alpha N}}{2(1-\beta^{-\alpha})}
 \quad\text{uniformly in }x.
\end{equation}
The recurrence requires $O(N)$ arithmetic operations per point. For $\beta=2$, one has $r_n=1$ and $y_n=\{2^nx\}$, so \eqref{eq:explicit} recovers the usual Takagi formula. The bound \eqref{eq:tailnum} controls truncation of the series; rounding errors in a numerical implementation require separate control.

\section{Spectra and classes of bases}\label{sec:spectra}
If the greedy expansion of the number $1$ is infinite, write $R_0=1$, $R_j=T_\beta^j(1)>0$,
\[
 \lambda(\beta)=\limsup_{j\to\infty}
       \frac{-\log R_j}{j\log\beta}.
\]
If the expansion is finite, set $\lambda(\beta)=0$. Here the map at $1$ is extended by the greedy formula; cylinders are still defined on $[0,1)$ and then closed. The parameter $\lambda(\beta)$ measures an additional exponential decrease in cylinder lengths.

\begin{lemma}\label{lem:geometry-general}
Suppose that the expansion of $1$ is infinite. Let $\EE_m$ be the set of right endpoints of full cylinders of order $m$, with $\EE_0=\{1\}$. Every nonfull cylinder of order $n$ has the form
\[
 I=[e-\beta^{-(m+j)}R_j,e],
\]
$e\in\EE_m$, $n=m+j$, $j\ge1$. Conversely, this interval is a cylinder of order $m+j$ for every
$e\in\EE_m$ and $j\ge1$. Moreover, $\#\EE_m\le\beta^m$,
\[
 \dist(x,\EE_m)\le(m+2)\beta^{-m}
 \quad(m\ge1,\ x\in[0,1]).
\]
\end{lemma}
\begin{proof}
If $I$ is the closure of a cylinder of order $n$, its half-open version has image $[0,\beta^n|I|)$ under $T_\beta^n$. Its subdivision consists of full branches with image $[0,1)$ and at most one incomplete branch, at the right endpoint. For a nonfull cylinder, consider its last full ancestor, of order $m$. The next $j$ steps follow incomplete branches and reproduce the first $j$ digits of the expansion of $1$. The cylinder shares the right endpoint of this ancestor and has length $\beta^{-(m+j)}R_j$. The converse follows by mapping the full cylinder affinely onto $[0,1)$. Disjoint interiors give $\#\EE_m\le\beta^m$. Finally, among any $m+1$ consecutive cylinders of order $m$, at least one is full (Bugeaud--Wang \cite{BW}, Theorem~1.2). Since the mesh is at most $\beta^{-m}$, the distance bound follows.
\end{proof}

\begin{theorem}\label{thm:general}
Let $\beta>1$, $0<\alpha\le1$ and $0\le\lambda=\lambda(\beta)<\infty$.
\begin{enumerate}
\item If $\lambda=0$, then $h_{\TT_{\beta,\alpha}}(x)=\alpha$ for every $x\in[0,1]$.
\item If $0<\lambda<\infty$, then
\[
 D_{\TT_{\beta,\alpha}}(u)=
 \begin{cases}
 \displaystyle\frac{\frac{(1+\lambda)}{\alpha}u-1}{\lambda}
       &\text{if }\displaystyle \frac{\alpha}{1+\lambda}\le u\le\alpha \\
 -\infty &\text{otherwise}.
 \end{cases}.
\]
For $\alpha(1+\lambda)^{-1}<u<\alpha$, put
\[
d=\lambda^{-1}(\frac{(1+\lambda)}{\alpha}u-1).
\]
Every nonempty open interval $U\subset(0,1)$ satisfies
\[
 \HH^d\bigl(U\cap\{x\in[0,1]:h_{\TT_{\beta,\alpha}}(x)=u\}\bigr)=\infty.
\]
The set
\[
\{x\in[0,1]:h_{\TT_{\beta,\alpha}}(x)=\alpha\}
\]
has full Lebesgue measure. The set
\[
\{x\in[0,1]:h_{\TT_{\beta,\alpha}}(x)=\frac{\alpha}{1+\lambda}\}
\]
is a dense $G_\delta$ of Hausdorff dimension zero.
\end{enumerate}
\end{theorem}

\begin{proof}
Write $h(x)=h_{\TT_{\beta,\alpha}}(x)$ and $L=1+\lambda$. If the expansion of the number $1$ is finite, cylinder lengths are uniformly comparable to $\beta^{-n}$: their images have length $1$ or one of the finitely many positive remainders. Corollary~\ref{thm:holder} then gives $h(x)=\alpha$ for every $x\in[0,1]$. Assume henceforth that the expansion is infinite.

We first determine the range of exponents. For every $\varepsilon>0$, there is $C(\beta,\varepsilon)\in(0,1]$ such that $R_j\ge C(\beta \varepsilon)\beta^{-(\lambda+\varepsilon)j}$ for all $j\ge1$. Lemma~\ref{lem:geometry-general} gives
\[
 |I|\ge C(\beta,\varepsilon)\beta^{-(L+\varepsilon)n}
 \qquad(I\in\CC_n(\beta)).
\]
It follows that $h(x)\ge\alpha(L+\varepsilon)^{-1}$ for every $x$, and hence $\alpha/L\le h(x)\le\alpha$. This proves the first assertion when $\lambda=0$. At $1$, the terminal cylinder of order $j$ has length $\beta^{-j}R_j$. Corollary~\ref{thm:holder} therefore gives $h(1)\le\alpha/L$, so $h(1)=\alpha/L$.

Suppose now that $\lambda>0$. To bound the dimensions of sublevel sets, fix $1<q<L$ and set
\[
 A_q=\{x\in[0,1]:\delta_n(x)<\beta^{-qn}
                    \text{ for infinitely many }n\}.
\]
A cylinder realizing this inequality cannot be full. In the notation of Lemma~\ref{lem:geometry-general},
\[
 C(\beta,\varepsilon)\beta^{-(m+(L+\varepsilon)j)}
 \le |I|<\beta^{-q(m+j)}.
\]
Thus
\[
 m\le c(\beta,q,\varepsilon)j+K(\beta,q,\varepsilon),\qquad
 c(\beta,q,\varepsilon)=\frac{L+\varepsilon-q}{q-1},\qquad
 K(\beta,q,\varepsilon)=
       \frac{\log_\beta(1/C(\beta,\varepsilon))}{q-1}.
\]
Also,
\[
|x-e|\le\dist(x,I)+|I|<\beta^{-q(m+j)}.
\]
The set $A_q$ is therefore contained in the limsup of the balls
\[
 B(e,\beta^{-q(m+j)}),\quad
 0\le m\le c(\beta,q,\varepsilon)j+K(\beta,q,\varepsilon),
 \quad e\in\EE_m.
\]
Here and below balls are intersected with $[0,1]$ when necessary. In this covering, $m+j\to\infty$ implies $j\to\infty$. The sum of the $s$th powers of the diameters is bounded by
\[
 2^s\sum_{j\ge1}
 \sum_{0\le m\le c(\beta,q,\varepsilon)j+K(\beta,q,\varepsilon)}
 \beta^{(1-qs)m-qsj}.
\]
This series converges if
\[
 s>\frac{c(\beta,q,\varepsilon)}
         {q(c(\beta,q,\varepsilon)+1)}
   =\frac{L+\varepsilon-q}{q(\lambda+\varepsilon)}.
\]
Letting $\varepsilon\downarrow0$ yields
\[
\dim_H A_q\le\frac{L-q}{q\lambda}.
\]
For $\alpha/L\le u<\alpha$, the pointwise characterization gives
\[
\{x\in[0,1]:h(x)\le u\}\subset A_q
\]
whenever $1<q<\alpha/u$. Letting $q\uparrow\alpha/u$, we obtain
\begin{equation}\label{eq:sublevels}
 \dim_H\{x\in[0,1]:h(x)\le u\}
       \le\frac{Lu/\alpha-1}{\lambda}.
\end{equation}

For the lower bound, fix $\alpha/L<u<\alpha$ and put
\[
 d=\frac{Lu/\alpha-1}{\lambda}\in(0,1),\qquad
 c=\frac{Ld}{1-d}.
\]
Then $d=c(c+L)^{-1}$ and $u/\alpha=(c+1)(c+L)^{-1}$. Choose a strictly increasing sequence $j_k$ such that $-\log_\beta R_{j_k}/j_k\to\lambda$, and set $m_k=\lfloor cj_k\rfloor$. Since $c>0$, we have $m_k\to\infty$. After deleting finitely many initial terms and reindexing, assume that $m_k\ge1$ for every $k\ge1$. Define
\[
 \ell_k=\beta^{-(m_k+j_k)}R_{j_k},\qquad
 r_k=\max\left\{
       \bigl(3(m_k+2)\beta^{-m_k}\bigr)^{1/d},\,\ell_k\right\}.
\]
Both terms in this maximum have the same exponential decay rate, so
\begin{equation}\label{eq:r-asympt}
 \frac{-\log_\beta r_k}{j_k}\to c+L.
\end{equation}
By Lemma~\ref{lem:geometry-general}, the balls centred at $\EE_{m_k}$ with radius $r_k^d$ cover $[0,1]$ for every $k$. The mass transference principle of Beresnevich--Velani \cite{BV}, Theorem~3 with $X=[0,1]$, $g(r)=r$ and $f(r)=r^d$, gives
\[
 Z_u=\limsup_k\bigcup_{e\in\EE_{m_k}}B(e,r_k),
 \qquad \HH^d(U\cap Z_u)=\infty
\]
for every nonempty open interval $U\subset(0,1)$. There are finitely many balls at each order and $r_k\to0$, so an enumeration of all the balls satisfies the hypotheses of that theorem. If $x\in B(e,r_k)$, the cylinder of order $m_k+j_k$ adjoining $e$ on the left has length $\ell_k\le r_k$. Hence $\delta_{m_k+j_k}(x)\le2r_k$. Equations~\eqref{eq:holder} and \eqref{eq:r-asympt} imply
\[
Z_u\subset\{x\in[0,1]:h(x)\le u\}.
\]
Now $\{x\in[0,1]:h(x)<u\}$ is the union of the sublevel sets $\{x\in[0,1]:h(x)\le v\}$ over rational $\alpha/L<v<u$. By \eqref{eq:sublevels}, each has dimension strictly less than $d$ and thus zero $\HH^d$ measure. Their union also has zero $\HH^d$ measure. Removing it leaves
\[
\HH^d(U\cap\{x\in[0,1]:h(x)=u\})=\infty,
\]
proving the claimed spectrum at every interior exponent.

It remains to consider the endpoints. The level $\alpha/L$ contains $1$ and has dimension zero by \eqref{eq:sublevels}. The level $\alpha$ has full Lebesgue measure, since its complement is a countable union of sublevel sets of dimension less than $1$. Every right endpoint $e$ of a full cylinder of order $m$ satisfies $\delta_{m+j_k}(e)\le\beta^{-(m+j_k)}R_{j_k}$, hence $h(e)=\alpha/L$. These endpoints are dense. Finally,
\[
 \{x\in[0,1]:h(x)=\alpha/L\}=
 \bigcap_{\substack{q\in\mathbb Q\\1<q<L}}
 \bigcap_{M\ge1}\ \bigcup_{n\ge M}
 \{x\in[0,1]:\delta_n(x)<\beta^{-qn}\}.
\]
Each union is relatively open and contains all right endpoints of full cylinders. It is therefore dense, proving the $G_\delta$ assertion.
\end{proof}

\label{sec:arithmetic}
A Parry base is a base for which the greedy expansion of $1$ is finite or eventually periodic. A Pisot number is an algebraic integer $\beta>1$ whose other conjugates all have modulus less than $1$. A Salem number is an algebraic integer $\beta>1$ whose other conjugates have modulus at most $1$, with at least one on the unit circle.

\begin{corollary}
\label{cor:parry}
If the greedy beta-expansion of $1$ is finite or eventually periodic, then, for $\alpha\in(0,1]$,
\[
 h_{\TT_{\beta,\alpha}}(x)=\alpha\quad(x\in[0,1]),
\]
and
\[
D_{\TT_{\beta,\alpha}}(u)=
\begin{cases}1&\text{if }u=\alpha\\
-\infty&\text{if }u\ne\alpha.
\end{cases}.
\]
\end{corollary}
\begin{proof}
In the infinite case, the remainders $R_n$ are the values of the successive suffixes of the expansion of $1$. Eventual periodicity implies that they take only finitely many values. All are positive, so $\inf_nR_n>0$ and $\lambda(\beta)=0$. In the finite case, our convention is $\lambda(\beta)=0$, and the cylinder lengths are uniformly comparable to $\beta^{-n}$, as above. Theorem~\ref{thm:general} applies in both cases.
\end{proof}

\begin{lemma}
\label{lem:salem-norm}
Let $\beta$ be a Salem number of degree $d$. There is a constant $c_\beta>0$ such that every nonzero remainder $R_n=T_\beta^n(1)$ satisfies
\begin{equation}\label{eq:salem-remainder}
 R_n\ge c_\beta(n+1)^{-(d-2)}.
\end{equation}
If $\beta$ is a Pisot number, its nonzero remainders are bounded below by a positive constant.
\end{lemma}
\begin{proof}
Let $(a_k)_{k\ge1}$ be the digits of the greedy beta-expansion of $1$, and set
\[
 P_n(X)=X^n-\sum_{k=1}^n a_kX^{n-k},
 \]
$R_n=P_n(\beta)$. The digits satisfy $0\le a_k\le A:=\lfloor\beta\rfloor$. The element $R_n$ belongs to the ring of integers of $K=\mathbb Q(\beta)$. If nonzero, its norm is a nonzero integer, and therefore
\begin{equation}\label{eq:remainder-norm}
 1\le |N_{K/\mathbb Q}(R_n)|
   =R_n\prod_{\substack{\gamma\text{ conjugate to }\beta\\\gamma\ne\beta}}
                        |P_n(\gamma)|.
\end{equation}
For a conjugate with $|\gamma|<1$,
\[
 |P_n(\gamma)|
 \le|\gamma|^n+A\sum_{r=0}^{n-1}|\gamma|^r
 \le1+\frac{A}{1-|\gamma|}.
\]
For $|\gamma|=1$, we have $|P_n(\gamma)|\le1+An$. A Salem number has the conjugate $\beta^{-1}$ and $d-2$ conjugates on the unit circle. This is Salem's structural result \cite{Salem,GH}. Thus \eqref{eq:remainder-norm} gives
\[
 1\le C_\beta R_n(1+An)^{d-2},
 \qquad C_\beta=1+\frac{A}{1-\beta^{-1}}.
\]
Taking $c_\beta=[C_\beta\max\{1,A\}^{d-2}]^{-1}$ proves \eqref{eq:salem-remainder}. For a Pisot number, all other conjugates have modulus less than $1$, so their contributions to the product are bounded independently of $n$. This gives a constant lower bound. For integer $\beta>1$, the expansion is finite and the empty product causes no difficulty.
\end{proof}

\begin{corollary}
\label{cor:salem}
If $\beta$ is a Salem or Pisot number, then, for $\alpha\in (0,1]$,
\[
 h_{\TT_{\beta,\alpha}}(x)=\alpha\quad(x\in[0,1]),
\]
and
\[
 D_{\TT_{\beta,\alpha}}(u)=
 \begin{cases}
 1&\text{if }u=\alpha\\
 -\infty&\text{if }u\ne\alpha
 \end{cases}.
\]
\end{corollary}
\begin{proof}
If the expansion of the number $1$ is finite, apply Corollary~\ref{cor:parry}.
Otherwise all remainders are positive and at most $1$. For a Salem
number, Lemma~\ref{lem:salem-norm} gives
\[
 0\le\frac{-\log R_n}{n\log\beta}
 \le\frac{(d-2)\log(n+1)-\log c_\beta}{n\log\beta}
 \to0.
\]
Thus $\lambda(\beta)=0$. The constant lower bound gives the same
conclusion for Pisot numbers. The result follows from
Theorem~\ref{thm:general}.
\end{proof}

\begin{remark}
Corollary~\ref{cor:salem} requires no eventual periodicity assumption on the expansion of $1$ for Salem numbers. In particular, it does not depend on whether the base is Parry.
\end{remark}

For a Salem number with an infinite expansion, the cylinder length formula and \eqref{eq:salem-remainder} imply
\[
 c_\beta(n+1)^{-(d-2)}\beta^{-n}
 \le |I|\le\beta^{-n}\qquad(I\in\CC_n(\beta)).
\]
Indeed, a nonfull cylinder has length $\beta^{-n}R_j$ with $1\le j\le n$, whereas a full cylinder has length $\beta^{-n}$. The estimate does not assert that the polynomial loss is attained. It rules out an exponential loss and gives
\[
 -\log\delta_n(x)=n\log\beta+O_\beta(\log(n+1))
 \quad\text{uniformly in }x.
\]
Thus a factor $C_n\to0$ in a length $C_n\beta^{-n}$ need not change a H\"older exponent: the condition $-\log C_n=o(n)$ leaves the leading logarithmic order unchanged. Exponential compression, measured here by $\lambda(\beta)>0$, produces the nontrivial spectra below.

\section{An explicit non-monofractal family}\label{sec:family}
For each integer $\ell\ge2$, the following proposition defines a base $\beta_\ell$. The associated function is the specialization of \eqref{eq:def} to this base:
\[
 \TT_{\beta_\ell,\alpha}(x)
 =\sum_{n\ge0}\beta_\ell^{-\alpha n}
       \sum_{I\in\CC_n(\beta_\ell)}\tau_I(x).
\]
We compute its spectrum and then use the same family to show why neighbouring cylinders are needed in the pointwise characterization.

\begin{proposition}\label{prop:family}
For every integer $\ell\ge2$, there is a unique $\beta_\ell\in(1,2)$ such that
\[
 \sum_{k=0}^{\infty}\beta_\ell^{-\ell^k}=1.
\]
The digits equal to $1$ in the greedy $\beta_\ell$-expansion of $1$ occur exactly at positions $\ell^k$. Moreover, $\lambda(\beta_\ell)=\ell-1$ and, for $0<\alpha\le1$,
\[
 D_{\TT_{\beta_\ell,\alpha}}(u)=
 \begin{cases}
 \displaystyle\frac{\frac{\ell}{\alpha}u-1}{\ell-1}&\text{if }\alpha/\ell\le u\le\alpha\\
 -\infty&\text{otherwise}
 \end{cases}.
\]
\end{proposition}
\begin{proof}
The sum is continuous and strictly decreasing on $(1,\infty)$, tends to infinity at $1$, and is less than $1$ at $2$. At its root $b=\beta_\ell$, one has $b^{-1}+b^{-\ell}<1$. In every proper suffix of the proposed digit sequence, successive $1$s are separated by at least $\ell$ positions. Its value is at most $b^{-1}(1-b^{-\ell})^{-1}<1$, which ensures that the expansion is greedy. For $N=\ell^k$ and $N\le j<\ell N$, the remainder is $R_j=b^j\rho_k$,
\[
 \rho_k=\sum_{r>k}b^{-\ell^r}\asymp b^{-\ell N}.
\]
Hence $-\log_bR_j=\ell N-j+O(1)$. Dividing by $j$ gives an upper limit of $\ell-1$, attained along $j=N$ as $N\to\infty$. Apply Theorem~\ref{thm:general}.
\end{proof}

For $\alpha=1$, the first examples are
\[
\begin{array}{c|c|c|c}
\ell&\beta_\ell\ \text{(approximate)}&\text{exponents}&D(u)\\ \hline
2&1.766398114550&[1/2,1]&2u-1\\
3&1.490890314609&[1/3,1]&(3u-1)/2\\
4&1.384450093664&[1/4,1]&(4u-1)/3
\end{array}
\]
These expansions are not eventually periodic, so the bases are not Parry. They are neither Pisot nor Salem, since $\lambda(\beta_\ell)=\ell-1>0$, whereas the proof of Corollary~\ref{cor:salem} gives $\lambda(\beta)=0$ for both classes.

The next example shows that $\delta_n(x)$ cannot be replaced in \eqref{eq:holder} by the length of the cylinder $I_n(x)$ containing $x$. These lengths can be comparable to $\beta^{-n}$ even when the pointwise exponent is strictly less than $\alpha$.

Fix $\ell\ge2$, let $b=\beta_\ell$, and set $M_j=\ell^{\ell^j}$,
\[
x_\infty=\sum_{j\ge1}b^{-M_j}.
\]
The $1$s are separated by at least $\ell-1$ zeros. The sequence is admissible, since every suffix has value at most $b^{-1}(1-b^{-\ell})^{-1}<1$. After a digit $1$ from a full image, successive zeros give image lengths $\min\{1,b^s(b-1)\}$ for $s\ge0$. Since $b^{\ell-1}(b-1)>1$, a full image is recovered before the next $1$. Consequently
\[
 (b-1)b^{-n}\le |I_n(x_\infty)|\le b^{-n}.
\]
Nevertheless, $h_{\TT_{b,\alpha}}(x_\infty)=\alpha/\ell$. Indeed, $e_j=\sum_{r\le j}b^{-M_r}$ is the right endpoint of a full cylinder of order $M_j$, by Lemma~\ref{lem:lower-digit} applied to the prefix of length $M_j$ of the expansion of $x_\infty$. For $N=M_{j+1}/\ell$, the cylinder of order $M_j+N$ adjoining $e_j$ on the left has length $O(b^{-(M_j+M_{j+1})})$. Since $x_\infty-e_j=O(b^{-M_{j+1}})$,
\[
 \delta_{M_j+M_{j+1}/\ell}(x_\infty)=O(b^{-M_{j+1}}).
\]
Corollary~\ref{thm:holder} and $M_j/M_{j+1}\to0$ give the upper bound $\alpha/\ell$; the reverse bound holds globally. The point $x_\infty$ is not a cylinder endpoint, since its greedy expansion is infinite.

\section{Comparison with the spectrum of Fan--Wang}\label{sec:comparison}
\label{sec:FW-audit}
Fan and Wang \cite{FW} determine the dimensions of the sets $\{x:Q(x)=q\}$, where $Q(x)$ is defined below from the lengths of the cylinders containing $x$. We examine what this result provides for the H\"older spectrum of $\TT_{\beta,\alpha}$. Our pointwise characterization involves $Q_*(x)$, defined from $\delta_n(x)$, so the relevant comparison is between the exact level sets of $Q$ and $Q_*$.

Fix $\beta>1$ and $\alpha\in(0,1]$, with $0<\lambda=\lambda(\beta)<\infty$. For $x\in[0,1)$, let $I_n(x)$ be the closure of the half-open cylinder specified by its first $n$ greedy digits, as in Section~\ref{sec:cylinders}. At an interior boundary this selects the interval to the right; for $x=1$, choose the last interval of the partition. Set
\[
 Q(x)=\limsup_{n\to\infty}\frac{-\log_\beta|I_n(x)|}{n},
 \qquad
 Q_*(x)=\limsup_{n\to\infty}\frac{-\log_\beta\delta_n(x)}{n}.
\]
The first quantity depends only on cylinders containing $x$. The second also depends on neighbouring cylinders: a short interval $I$ near $x$ can satisfy $|I|+\dist(x,I)<|I_n(x)|$ and contribute to a larger value of $Q_*(x)$. Since $\delta_n(x)\le|I_n(x)|$, we have $Q_*(x)\ge Q(x)$ for all $x\in[0,1]$.

Our parameter $\lambda(\beta)$ agrees with that of Fan--Wang. To see this, let $(a_i)_{i\ge1}$ be the infinite expansion of the number $1$ and set
\[
 t_n=\min\{k\ge0:a_{n+k+1}>0\}\qquad(n\ge1).
\]
Thus $t_n$ counts the zeros immediately following $a_n$. Since $R_n=\sum_{i\ge1}a_{n+i}\beta^{-i}$,
\[
 \beta^{-(t_n+1)}\le R_n
 \le\frac{\lfloor\beta\rfloor}{\beta-1}\,\beta^{-t_n}.
\]
Consequently $-\log_\beta R_n=t_n+O_\beta(1)$ and
\[
\lambda(\beta)=\limsup_{n\to\infty}\frac{t_n}{n}.
\]
Lemma~6.2 of \cite{FW} identifies this upper limit with their parameter $\lambda(\beta)$ when it is positive, as assumed here.

To avoid endpoint coding conventions, let $\mathcal B_\beta$ be the set of all cylinder endpoints and put $\Omega_\beta=[0,1]\setminus\mathcal B_\beta$. The set $\mathcal B_\beta$ is countable. On $\Omega_\beta$, the containing cylinder is unique, and $Q(x)$ is the exponent studied in \cite{FW}. Their Theorem~1.3 gives, for $1<q<1+\lambda$,
\[
 S(q):=\dim_H\{x\in\Omega_\beta:Q(x)=q\}
      =\frac{1+\lambda-q}{q\lambda}.
\]
These dimensions are positive and hence unchanged by removal of the countable endpoint set. Thus $S$ denotes the dimension function for the level sets of $Q$, rather than the H\"older spectrum of our function.

Our characterization instead gives
\[
 h_{\TT_{\beta,\alpha}}(x)=\frac{\alpha}{Q_*(x)}.
\]
For $u>0$, the H\"older spectrum is therefore the dimension of $\{x:Q_*(x)=\alpha/u\}$. The inequality $Q_*(x)\ge Q(x)$ does not identify this set with an exact level set of $Q$. The following example locates the difficulty within Fan--Wang's construction itself.

Their construction gives a set $E$ whose expansions alternate prescribed blocks and free blocks chosen from a family of admissible words. Every $x\in E$ satisfies $Q(x)=q$, and the number of available free blocks gives a lower bound for $\dim_HE$. Using this set directly for our function would require $Q_*(x)=q$ throughout $E$. This need not hold.

\begin{proposition}\label{prop:FW-moran}
For the base $b=\beta_2$ defined by $\sum_{r\ge0}b^{-2^r}=1$, the construction in Section~6.4 of \cite{FW}, at level $q=6/5$, can be chosen to contain a point $x$ such that
\[
 Q(x)=\frac65,\qquad Q_*(x)=2,\qquad
 h_{\TT_{b,\alpha}}(x)=\frac{\alpha}{2}.
\]
Its Moran set is therefore not necessarily contained in $\{x\in[0,1]:Q_*(x)=6/5\}$.
\end{proposition}
\begin{proof}
We first construct a point allowed by Fan--Wang's blocks, which gives $Q(x)=6/5$, and then find neighbouring cylinders forcing $Q_*(x)=2$. Write $(a_i)_{i\ge1}$ for the expansion of $1$ in base $b$: $a_i=1$ precisely when $i$ is a power of two. For $m=2^r$, the zero block following $a_m$ has length $t_m=m-1$, and
\[
 R_m\asymp b^{-m}.
\]
Fix the integers $N,M$ in Fan--Wang's construction. Their family $F$ of free words contains $0^M$: the zero word is admissible in their subsystem and satisfies the terminal-zero condition.

Choose powers of two $m_k$ such that $m_k/k\ge\sum_{j<k}m_j$ and $m_{k-1}/m_k\to0$. The positions $h_k$ can be chosen so that $h_k=4m_k+O(M)$,
\[
h_k-(h_{k-1}+2m_{k-1})\in M\mathbb N\quad(k\ge2),
\]
with $h_1\in M\mathbb N$. The gaps are positive if the $m_k$ are sufficiently far apart. These choices satisfy conditions~(6.3)--(6.4) of \cite{FW}, since
\[
 \frac{h_k+m_k+t_{m_k}}{h_k+m_k}\to\frac65.
\]
The prescribed block $G_k=(a_1,\ldots,a_{m_k},0^{m_k})$ has length $2m_k$ and starts at position $h_k+1$.
Choose $0^M$ at every free position, and let $x$ be the resulting point. It belongs to $E$, and Lemma~6.3 of \cite{FW} gives $Q(x)=6/5$. The resulting greedy expansion has infinitely many nonzero digits and arbitrarily long zero blocks. In particular, $x$ is not an endpoint, and the coding conventions agree there.

Set $p_k=h_k+m_k$ and
\[
 e_k=\sum_{i=1}^{p_k}a_i(x)b^{-i}.
\]
The digit at position $p_k$ is $1$. All subsequent digits are zero up to position $h_{k+1}+1$, where the first $1$ of the next block occurs. Thus
\[
 b^{-(h_{k+1}+1)}\le x-e_k
 \le\frac{b^{-h_{k+1}}}{b-1},
 \qquad x-e_k\asymp b^{-h_{k+1}}.
\]
Lemma~\ref{lem:lower-digit}, applied to the prefix of length $p_k$, shows that $e_k$ is the right endpoint of a full cylinder of order $p_k$.

Within this full cylinder, the terminal cylinder of depth $m_{k+1}$ is an interval $J_k$ adjoining $e_k$ on the left. Its order is $N_k=p_k+m_{k+1}$ and its length is
\[
 |J_k|=b^{-(p_k+m_{k+1})}R_{m_{k+1}}
       \asymp b^{-(p_k+2m_{k+1})}.
\]
Since $x>e_k$,
\[
 \delta_{N_k}(x)\le x-e_k+|J_k|
 \lesssim b^{-h_{k+1}}+b^{-(p_k+2m_{k+1})}.
\]
Now $h_{k+1}\sim4m_{k+1}$ and $p_k=o(m_{k+1})$, so the first term is negligible compared with the second. Hence
\[
 Q_*(x)\ge
 \lim_{k\to\infty}\frac{p_k+2m_{k+1}+O(1)}{p_k+m_{k+1}}=2.
\]
The global bound $Q_*(x)\le1+\lambda(b)=2$ gives equality. The pointwise characterization then yields $h_{\TT_{b,\alpha}}(x)=\alpha/2$.
\end{proof}

Proposition~\ref{prop:FW-moran} rules out identifying $Q(x)$ and $Q_*(x)$ pointwise on Fan--Wang's set. It does not assert that their exact level sets have different dimensions.

For $1<q<1+\lambda$, the inequality $Q_*(x)\ge Q(x)$ and Fan--Wang's result give
\[
 \dim_H\{x\in[0,1]:Q_*(x)\ge q\}\ge S(q).
\]
To obtain the exact level $Q_*(x)=q$, one must remove points with $Q_*(x)>q$. Dimension bounds alone do not justify this step, since $\sup_{q<r<1+\lambda}S(r)=S(q)$. In the proof of Theorem~\ref{thm:general}, mass transference gives, for $d=S(q)$, a set
\[
Z\subset\{x\in[0,1]:Q_*(x)\ge q\}
\]
with $\HH^d(Z\cap U)=\infty$ in every nonempty open interval $U\subset(0,1)$. For each rational number $r\in(q,1+\lambda)$, the upper bound gives
\[
 \dim_H\{x\in[0,1]:Q_*(x)\ge r\}\le S(r)<d.
\]
Since $Q_*(x)\le1+\lambda$, the set $\{x\in[0,1]:Q_*(x)>q\}$ is the countable union of these sets and has zero $\HH^d$ measure. Removing it preserves the infinite $\HH^d$ measure in every $U$, which gives the required lower bound for the exact level $Q_*(x)=q$.


\begin{thebibliography}{99}
\bibitem{Allaart}
P. Allaart,
\emph{The pointwise H\"older spectrum of general self-affine functions on an interval},
\href{https://arxiv.org/abs/1907.09660v2}{arXiv:1907.09660v2} (2019).

\bibitem{AJ}
A. Ayache and S. Jaffard,
\emph{H\"older exponents of arbitrary functions},
Revista Matemática Iberoamericana 26 (2010), no.~1, 77--89.

\bibitem{BDK}
E. Bayraktar, P. Das and D. Kim,
\emph{H\"older regularity and roughness: construction and examples},
Bernoulli 31 (2025), 1084--1113.

\bibitem{BS}
M. Ben Slimane, M. Ben Abid, I. Ben Omrane and B. Halouani,
\emph{Criteria of pointwise and uniform directional Lipschitz regularities on tensor products of Schauder functions},
Journal of Mathematical Analysis and Applications 460 (2018), 496--515.

\bibitem{BV}
V. Beresnevich and S. Velani,
\emph{A Mass Transference Principle and the Duffin--Schaeffer conjecture for Hausdorff measures},
Annals of Mathematics 164 (2006), 971--992.

\bibitem{BW}
Y. Bugeaud and B.-W. Wang,
\emph{Distribution of full cylinders and the Diophantine properties of the orbits in $\beta$-expansions},
Journal of Fractal Geometry 1 (2014), no.~2, 221--241.

\bibitem{CN}
M. Clausel and S. Nicolay,
\emph{Some prevalent results about strongly monoH\"older functions},
Nonlinearity 23 (2010), 2101--2116.

\bibitem{FW}
A.-H. Fan and B.-W. Wang,
\emph{On the lengths of basic intervals in beta expansions},
Nonlinearity 25 (2012), 1329--1343.

\bibitem{JN}
S. Jaffard and S. Nicolay,
\emph{Pointwise smoothness of space-filling functions},
Applied and Computational Harmonic Analysis 26 (2009), no.~2, 181--199.

\bibitem{JS}
J. Jaerisch and H. Sumi,
\emph{Spectral gap property for random dynamics on the real line and multifractal analysis of generalised Takagi functions},
Communications in Mathematical Physics 377 (2020), 1--36.

\bibitem{GH}
E. Ghate and E. Hironaka,
\emph{The arithmetic and geometry of Salem numbers},
Bulletin of the American Mathematical Society 38 (2001), no.~3, 293--314.

\bibitem{Salem}
R. Salem,
\emph{Algebraic Numbers and Fourier Analysis},
D. C. Heath and Company, Boston, 1963, p.~26.

\bibitem{Suzuki}
S. Suzuki,
\emph{A generalization of the Takagi function for beta-expansions},
preprint (2026).
\end{thebibliography}
\end{document}